\def\R{I \kern-.33em R}
\def\C{I \kern-.66em  C}
\def\fs{\footnotesize}

\def\al{\alpha}
\def\bt{\beta}
\def\ga{\gamma}
\def\Ga{\Gamma}

\def\ep{\epsilon}

\def\si{\sigma}

\def\pt{\partial}
\def\na{\nabla}

\def\uno{\hbox{\bf I}}

\def\h{\tilde h}
\def\g{\tilde g}
\def\J{J_o}
\def\zz{h^{-1} g}
\def\ba{\begin{array}}
\def\ea{\end{array}}

\documentstyle[11pt]{article}
\begin{document}

\title{ Almost	Complex and Almost  Product\\ Einstein Manifolds\\
from a Variational Principle}

\author{
Andrzej Borowiec\thanks{On leave  from the
Institute of Theoretical Physics, University of Wroc{\l}aw,
pl. Maksa Borna 9, 50-204  WROC{\L}AW (POLAND), E-mail:
borow@ift.uni.wroc.pl}\\
 Marco Ferraris\\
Mauro Francaviglia\\
Igor Volovich\thanks
{Permanent address: Steklov Mathematical
Institute, Russian Academy of Sciences,
Vavilov St. 42, GSP--1, 117966
MOSCOW (RUSSIA).}\\
{\fs Istituto di Fisica Matematica ``J.--L. Lagrange''}\\
{\fs Universit\`a di Torino}\\ {\fs Via C. Alberto 10, 10123 TORINO (ITALY)}
}
\date{\bf Preprint TO--JLL--P 7/96}
\maketitle

\begin{abstract}
It is shown that  the first order (Palatini) variational principle for
a generic nonlinear metric-affine Lagrangian depending
on the (symmetrized) Ricci square invariant leads
to an almost-product Einstein structure or to an almost-complex
anti-Hermitian Einstein structure on a manifold. It is proved that a real
anti-Hermitian metric on a complex manifold satisfies the K\"ahler condition
on the same manifold treated as a real manifold if and
only if the metric is the real part of a holomorphic metric.
A characterisation of anti-K\"ahler Einstein manifolds
and almost-product Einstein manifolds is obtained. Examples of such manifolds
are considered.
\end{abstract}
\newpage

\section{Introduction}

$~~~~$Almost-complex and almost-product
structures are among the most important geometrical structures
which can be considered
on a manifold \cite{Che,Yan,Kob}. The aim of this paper is to show that
 structures of this kind appear in a natural way from a variational principle
based on  a general class of Lagrangians depending on the Ricci
square invariant constructed out of a metric and a symmetric connection;
in particular we will show that an anti-Hermitian metric and its special
case, so called "anti-K\"ahlerian" metric, appear naturally from our
variational principle. Manifolds with such metrics are much less studied
than the familiar Hermitian and K\"ahlerian cases. We hope
that our variational principle and "the universality of the Einstein
equations" (see below) will provide an additional motivation for
investigating these manifolds.

Let $M$ be a  differentiable manifold of dimension $n$ and $L(M)$
be its frame bundle, a principal fibre bundle over $M$ with group $GL(n;\R)$.
Let $G$ be a Lie subgroup of $GL(n;\R)$. A differentiable subbundle $Q$
of $L(M)$ with structure group $G$ is called a {\it $G$-structure on $M$ }
\cite{Che,Kob}. The classification and integrability of $G$-structures
have been studied in differential geometry; algebraic-topological conditions
on $M$ which are necessary for the existence
of a $G$-structure on $M$ can be given in terms of characteristic
classes (see, for example, \cite{Sun}). We also recall that there
is a natural one-to-one correspondence between	pseudo-Riemannian
metrics of signature $q$ on $M$ and  $O(p,q;\R)$-structures on $M$, with
$p+q=n$. An $O(p,q;\R)$-structure  is integrable if and only if the
corresponding pseudo-Riemannian metric has vanishing Riemann curvature.
 If $G=GL(p;\R)\times GL(q;\R)$ then the $G$-structure
 is called an {\it almost-product structure} \cite{Yan,Rei,YK};
 if $G=O(r,s;\R)\times O(k,l;\R)$ then the $G$-structure is called
a {\it (pseudo-)Riemannian almost-product structure} \cite {Yan,Wal,Gra}.
If $n$ is even, $ n=2m$, and one considers $GL(m;\C)$ as a subgroup
of $GL(2m;\R)$ then a $GL(m;\C)$-structure is called an {\it
almost-complex structure} \cite {Yan,Kob}; if, moreover, one considers
$O(m;\C)$ as a subgroup of $GL(m;\C)$ then an $O(m;\C)$-structure
defines an {\it almost-complex anti-Hermitian structure}
\cite{Nor,GB,CHR,GI,BB}.

We will here use an equivalent description of these $G$-structures.
Let $M$ be a manifold and $P$ be an endomorphism of the tangent
bundle $TM$ satisfying $P^2=I$, where $I=identity$. Then $P$
defines an almost-product structure on $M$.
If $g$ is a  metric on $M$ such that
 $g(PX, PY)= g(X, Y)$
for arbitrary vectorfields $X$ and $Y$ on $M$,
then the triple $(M,g,P)$ defines a
 (pseudo-) Riemannian almost-product structure.
 Geometric properties of (pseudo-) Riemannian almost-product structures have
been studied in \cite{Yan,Wal,Gra,Nav,GN,CM,Roc,Step,Rei}.
If, moreover, $g$ is an Einstein
 metric ({\it i.e.} $Ric(g)=\ga g$ holds, where $Ric(g)$ is the Ricci
 tensor and $\ga$ is a constant) then the triple
$ (M,g,P)$ shall be called an {\it almost-product Einstein manifold}.

Analogously, if $J$ is an endomorphism of the tangent
bundle	$TM$ satisfying $J^2=-I$, then $J$ defines
an almost-complex structure on $M$. An almost-complex structure is
integrable if and only if it comes from a complex structure (see \cite{KN}).
If $g$ is a  metric on $M$ such that $g(JX, JY)=- g(X, Y)$
for arbitrary vectorfields $X$ and $Y$ on $M$
then the triple $(M,g,J)$ defines   an almost-complex
anti-Hermitian structure. The metric $g$ in this case is called a
{\it Norden metric} and in complex coordinates it has the form
$$
ds^2=g_{ab}dz^adz^b+g_{\bar a\bar b}dz^{\bar a}dz^{\bar b}
$$
where  $g_{\bar a\bar b}={\bar g_{ab}}$. This canonical form
differs from the well known form of an Hermitian
metric $ds^2=2g_{a\bar b}dz^{ a}dz^{\bar b}$.
We will show (theorem 4.2) that the condition $\nabla J=0$,
 where $\nabla$ is the Levi-Civita connection, is equivalent in this case to
analyticity of the metric:$~$ ${\bar\partial} _{ c}g_{ab}=0$ . Such
anti-Hermitian metrics shall be called {\it anti-K\"ahlerian} metrics since
for an Hermitian metric the condition $\nabla J=0$ defines a K\"ahlerian
metric.

If $g$ is an Einstein metric, {\it i.e.} $Ric(g)=\ga g$ holds, then the
almost-complex anti-Hermitian manifold $(M,g,J)$ is called an
{\it anti-Hermitian Einstein manifold}. We will consider
an important particular class of such manifolds, namely those characterized by
{\it anti-K\"ahlerian Einstein metrics}.
Let us stress that we treat {\it the whole} complex manifold as a real
manifold and in this way we get a real
Einstein metric with signature $(m,m)$. Another approach to complex Einstein
equations, dealing with a {\it	real section} of a complex manifold and
aiming to get the Lorentz signature, has been considered in
\cite{Pen,Man,Bru,Hit,Gin,Fla}.

 These $G$-structures can be conveniently defined as a triple $(M,g,K)$,
 where $g$ is a metric on $M$ and $K$ is a (1,1) tensor
 field on $M$ such that $K^2=\ep I$ and $g(KX,KY)=\ep g(X,Y)$
  for arbitrary vectorfields $X$ and $Y$ on $M$ ($\ep\neq 0 $ is a real
constant).  We shall also call them {\it K-structures}.
If $\epsilon =1$ then $K$ defines
an almost-product structure on	$M$; if $\epsilon =-1$ then $K$ defines
an almost-complex structure on $M$.
The more general case $\ep >0$ can be reduced to $\ep =1$ by  a suitable
rescaling, while the case $\ep <0$ is reduced to $\ep =-1$.
 In any coordinate system one has
$K^{\mu}_{\alpha}K^{\alpha}_{\nu}=\epsilon\delta^{\mu}_{\nu}$
and  $K^t gK=\ep g$, where  $K^t $ is the transpose  matrix,
 $\mu ,\nu ,\alpha =1,2,...,n=dim\ (M)$ and $\delta^{\mu}_{\nu}$ is the
Kronecker symbol.

One can then define a new metric $h$ by the relation $h(X,Y)=g(KX,Y)$, or
equivalently $h=gK$, {\it i.e.} in local coordinates
$h_{\mu\nu}=g_{\mu\alpha}K^{\alpha}_{\nu}$; then the following holds
$$
(g^{-1}h)^2 =\ep I. \eqno	    (1.1)
$$
The relation (1.1) for $\ep =+1$ or $-1$ is equivalent
to $(h^{-1}g)^2=\ep I$ and there is a one-to-one correspondence
between the $G$-structure $(M,g,K)$ and the $G$-structure
$(M,h,K^{-1})$. Thence the $G$-structure given by the triple  $(M,g,K)$
can be equivalently described by the triple $(M,g,h)$,
where $g$ and $h$ are metrics on $M$ satisfying  (1.1).
We call such metrics {\it twin metrics} or {\it dual metrics}.

In this paper, starting from a manifold $M$ endowed  with a metric
$h=(h_{\mu\nu})$ and a symmetric linear connection
$\Ga =(\Ga_{\mu\nu}^{\alpha})$,
we obtain  a $K$-structure as
$$
K^{\al}_{\mu}=h^{\al\nu}S_{\mu\nu}
$$
where $S_{\mu\nu}\equiv R_{(\mu\nu)}(\Ga)$ is the symmetric part of the
Ricci tensor of the given connection $\Ga$. The general idea is the following.
Let us first set $g_{\mu\nu}=R_{(\mu\nu)}(\Ga)$; according to the results of
our earlier paper \cite{BFFV} one can
show that $g$ is in fact a new metric and that $h$ and $g$ are "twin metrics"
if one assumes a suitable variational principle based on the action
$$
A(\Gamma, h)=\int_M f(S)\sqrt h \ dx \eqno	    (1.2)
$$
and imposes independent variations over the metric and the connection.
Here $f(S)$ is a given function
of one real variable, which we assume to be  analytic, while
the scalar $S=S(\Gamma, h)$ is the Ricci square invariant
$$
S=h^{\mu \al}h^{\nu \bt}S_{\al \bt}S_{\mu \nu}
$$
If $f$ is a generic analytic function and $n>2$ one gets either a (pseudo-)
Riemannian almost-product structure
or an almost-complex structure. In fact, as we shall see below (Theorem 2.4),
 the Euler-Lagrange equations for (1.2) are generically equivalent to
 the following system of equations for two metrics
 $h_{\mu\nu}$ and
$g_{\mu\nu}$:
$$
(h^{-1}g)^2 =\frac{c}{n}I  \eqno (1.3)
$$
$$
Ric(g)=g      \eqno (1.4)
$$
where the real number $c$ is a root of the equation
$$
f^{\prime} (S)S-{n\over4}f(S)=0 \eqno (1.5)
$$
As an example, if one takes
$f(S)=(nS+c(8-n))^2$, then $S=c$ is a solution of (1.5) if $n\neq 8$ and
another (degenerate) solution is $S=c(n-8)/n$.
Turning to the general discussion, if $c>0$ then, as it was explained above,
solutions of (1.3) - (1.4) are in one-to-one correspondence
with almost-product Einstein manifolds $(M,g,P)$; while for
$c<0$ one gets	anti-Hermitian Einstein manifolds $(M,g,J)$.

Before proceeding further let us explain why the action (1.2) is interesting
and important in Mathematical Physics and especially in the theory of gravity.

As is well known  gravitational Lagrangians which are nonlinear in the scalar
curvature of a metric give rise to equations with higher (more than second)
derivatives or to the appearance of additional matter fields
\cite{Stelle,MFF}.
This strongly depends on having taken
a metric as the only basic variable and the equations ensuing
from such Lagrangians show an explicit dependence on
the Lagrangian itself. An important example of a non-linear
Lagrangian leading to equations with higher derivatives
is given by Calabi's variational principle \cite {Cal}, which shall be
discussed in a forthcoming paper.

It was	shown in \cite{FFV1} that, in contrast,  working in the first
order (Palatini) formalism, i.e. assuming  independent variations with respect
to a metric and a symmetric  connection, then, for a large
class of Lagrangians of the form $f(R)$,
where $R$ is the scalar curvature, the equations obtained are almost
independent on the	Lagrangian, the only such dependence being in fact
encoded into constants (cosmological and Newton's ones). In this sense the
equations obtained are ``universal'' and turn out to be Einstein equations in
generic cases. Considering nonlinear gravitational Lagrangians which  still
generate Einstein equations is particularly important since they provide a
simple but general approach to	governing
topology in dimension two \cite{FFV2} and in view of applications to string
theory \cite{Vol}.

In a previous paper of ours \cite{BFFV2} this discussion was
extended to the case of Lagrangians with an arbitrary dependence
on the square of the symmetrized Ricci tensor  of a metric and a
(torsionless) connection,
finding roughly  that the universality of Einstein  equations
also extends to this class of spacetimes. In this case, however,
new  important properties
appear: as we have already mentioned above, depending in fact on  the form of
the Lagrangian and on
the signature of the metric, one gets an almost-product Einsteinian
structure or an almost-complex Einsteinian structure. Topological and
geometrical obstructions for the global  existence of a solution of the
variational problem for this class of Lagrangians  will be here considered
in sections 4 and 5.

 Recently there has been
some interest on the problem of  signature change in General Relativity
\cite{Ell,Hay,DT}; the non-standard signature $(10+2)$ has been considered
also in superstring theory ($F$-theory
\cite{Vaf}) and extra time-like dimensions in  Kaluza-Klein
theory are considered in \cite{Sak,AV,ADV}.
Our results seem therefore to show new aspects of this problem which
can be relevant also for quantum gravity. For a mathematical
consideration of metrics with arbitrary signature
which can be relevant to mathematical physics see for example
\cite{Cer2,GS,Law,BBR}

This paper is organized as follows.
In the next section it will be shortly recalled how to get eqs. (1.3) and
(1.4) from the action (1.2).  In section 3 it is shown that eq. (1.3)
can be always solved locally for any  given metric $g$, in particular
satisfying (1.4). In section 4 we discuss the $K$-structures.
In  section 5 we discuss  the problems of the global existence as well as
the classification of almost-product Einstein manifolds.
In section 6 we prove that a real anti-Hermitian metric on a complex manifold
satisfies the K\"ahler condition on the same manifold treated as a real
manifold if and only if the metric is the real part of a holomorphic metric
on this manifold. Finally, we consider also examples of  almost-product
Einstein manifolds and anti-K\"ahlerian
Einstein manifolds. Theorems of section 3 are proved in the Appendix.

\section{ Field Equations}

 $~~~~$In this section we shall present in a more geometrical form the
results of our earlier paper \cite{BFFV2}, which form the basis of the further
 results presented hereafter.
 Consider,  in the first order (Palatini) formalism, the family of actions
$$
A(\Gamma, h)=\int_M f(S)\sqrt h \ dx \eqno	    (2.1)
$$
\noindent where: $M$ is a n-dimensional manifold $(n>2)$ endowed with a metric
$h_{\mu \nu}$ and a torsionless (i.e., symmetric) connection $\Ga_{\mu \nu}
^{\sigma}$; the Lagrangian density is $L=f(S)\sqrt h$,
where $f(S)$ is a given function
of one real variable, which we assume to be  analytic
and $\sqrt h$ is a shorthand for $|\det (h_{\mu\nu})|^{1\over 2}$;
the scalar $S$ is the symmetric part of Ricci square-invariant,
considered as a first order scalar concomitant of a metric and
(torsionless) connection, {\it i.e.}
$$
S = S(h, \Ga)=
h^{\mu \al}h^{\nu \bt}S_{\al \bt}S_{\mu \nu}	\eqno (2.2)
$$
being $S_{\mu \nu}=R_{(\mu \nu)}(\Ga)$	the symmetric part of Ricci tensor,
defined according to
$$
R_{\mu \nu \sigma}^{\lambda}(\Ga)
=\partial_{\nu} \Gamma_{\mu
\sigma}^{\lambda} - \partial_{\sigma}\Gamma_{\mu
 \nu}^{\lambda} + \Gamma_{\alpha\nu}^{\lambda} \Gamma_{\mu
\sigma}^{\alpha} - \Gamma_{\alpha \sigma}^{\lambda}
\Gamma_{\mu \nu}^{\alpha}
$$
$$
R_{\mu \sigma}(\Ga)=R_{\mu \nu \sigma}^{\nu}(\Ga) \ \ \ \
(\alpha, \mu, \nu,...=1,...,n)
$$

Following \cite{BFFV2} the Euler--Lagrange equations of the
action (2.1) with respect to
independent variations of
$h$ and $\Gamma$ are
$$
f^{\prime} (S)h^{\al \bt}S_{\mu \al}S_{\nu \bt}-{1\over 4}f(S)h_{\mu
\nu}=0 \eqno (2.3)
$$

$$
 \nabla_{\lambda}(f^{\prime} (S)\sqrt h \
h^{\mu \al}h^{\nu \bt}S_{\al \bt})=0	\eqno (2.4)
$$
\noindent where $ \nabla_{\lambda}$ is the covariant derivative with respect
to $\Gamma$.
Transvecting (2.4) with $h^{\mu \nu}$ tells  us that the scalar $S$  has to
obey the following real analytic equation
$$f^{\prime} (S)S-{n\over4}f(S)=0 \eqno (2.5)
$$
which allows to describe the general features of the non-linear system
(2.3)-(2.4) and tells us, in turn, that $S$ is generically forced to be a
constant. More precisely, it was shown in \cite{BFFV2} that whenever (2.5)
admits an (isolated) simple root $S=c$ then
the system is "essentially equivalent" to Einstein equations for a new metric
$g$ with a cosmological constant, in the precise sense which is hereafter
described  in greater detail. Consider in fact any solution
$$
S=c \eqno (2.6)
$$
\noindent of eq. (2.5) and assume that
$f^{\prime}(c) \neq 0$. Then eq. (2.4) reduces to
$$
\nabla_{\lambda}(\sqrt h h^{\mu \al}h^{\nu \bt}S_{\al \bt})=0 \eqno (2.7)
$$
\noindent while equation (2.3) reduces to
$$
h^{\al \bt}S_{\mu \al}S_{\nu \bt}=\ep h_{\mu \nu}
\eqno (2.8)
$$
\noindent where a new constant $\ep $ depending on $c$ arises according to
the rule
$$
\ep  =f(c)/4f^{\prime}(c)=c/n  \eqno (2.9)
$$
From (2.8) the regularity condition
$$
[\det(S_{\mu \nu})]^2=\ep^n[\det(h_{\mu \nu})]^2
\eqno (2.10)
$$
follows, which entails in particular that $det(S_{\mu \nu}) \neq 0$ provided
$\ep\neq 0$.
Under this last hypothesis let $S^{\mu \nu}$ be the inverse matrix of
$S_{\mu \nu}$ , so that from (2.8) we have:
$$
h^{\mu\al }h^{\nu\bt}S_{\al\bt}=\ep S^{\mu \nu}
\eqno (2.11)
$$
By using (2.10) and (2.11) we finally rewrite (2.7) as follows
$$
\nabla_{\lambda}[\sqrt {|det (S_{\al\bt}(\Ga))|}S^{\mu\nu}(\Ga)]
 =0 \eqno(2.12) $$
which will be considered  as a new equation in $\Ga$.

Let us recall now the following well-known result, essentially due to
Levi-Civita: for $n>2$, any metric $g$ and any symmetric connection $\Ga$
the general solution of the equation
$$
\na_{\al}(\sqrt g\,g^{\mu \nu})=0 \eqno (2.13)
$$
considered as an equation for $\Ga$ is	the Levi--Civita connection
$\Ga=\Ga_{LC}(g)$ {\it i.e.}
locally $$
\Ga^\si_{\mu \nu}(g)= {1\over2}g^{\si \al}(\pt_\mu g_{\nu \al} +
\pt_\nu g_{\mu \al} - \pt_\al g_{\mu \nu})\eqno (2.14)
$$
Therefore, the Ricci tensor $R_{\mu \nu}(\Ga)$ of $\Ga$ is
automatically symmetric and in fact
identical to the Ricci tensor $R_{\mu \nu}(g)$ of the metric $g$
itself.\smallskip\\
We can then prove the following:\smallskip\\

{\bf Proposition 2.1.}\ {\em Let us assume that $det (S_{\al\bt})\neq 0$.
Then a connection $\Ga$ satisfies  eq.
(2.12) if and only if there exists a metric $g_{\mu \nu}$ such that
$$
R_{\mu\nu}(g) = g_{\mu\nu} \
\eqno (2.15)
$$
and $\Ga=\Ga_{LC}(g)$ is  the Levi--Civita connection of $g$}.

{\it Proof.} Let $\Ga$ be a connection	satisfying  eq. (2.12) and let us set
$$
 g_{\mu\nu} =S_{\mu\nu}(\Ga)
\eqno (2.16)
$$
The tensorfield $g$ is a metric due to the condition $det (S_{\al\bt})\neq 0$.
Then it follows that  $\Ga$ has to be the Levi--Civita connection
of the metric $g$; moreover one has
$$
 S_{\mu\nu}(\Ga)=R_{(\mu\nu)}(\Ga)=R_{\mu\nu}(g)
\eqno (2.17)
$$
so that (2.15) follows from (2.16) and (2.17).

Conversely, let us  give a metric $ g_{\mu\nu}$ satisfying  (2.15) and let us
take $\Ga=\Ga_{LC}(g)$. One has again relations (2.17). From (2.17) and
(2.15) it follows then that  (2.16) holds. Therefore
$ g^{\mu\nu} =S^{\mu\nu}(\Ga)$ and $\det (S_{\mu\nu}(\Ga))=
\det ( g_{\mu\nu})\neq 0$.
Hence we see that eq. (2.14) reduces to (2.15), which is satisfied
since  $\Ga=\Ga_{LC}(g)$. Our claim is then proved.\hfill (Q.E.D.)\\

 According to the previous discussion, we see
 that the Euler--Lagrange equations (2.3) and (2.4) are thence equivalent
to the following equations for two metrics $h$ and $g$:
$$
h^{\al \bt}g_{\mu \al}g_{\nu \bt}=\ep h_{\mu \nu}\eqno (2.18)
$$
$$
R_{\mu\nu}(g) = g_{\mu\nu} \  \
\eqno (2.19)
$$
which  are in fact nothing but eqs. (1.3) and (1.4) of the Introduction.
The relation between the system (2.18)-(2.19) and the Euler-Lagrange
equations in the form (2.3)-(2.4) or (2.11)-(2.12) is given by setting
$g_{\mu\nu}=R_{(\mu\nu)}(\Ga)$.

We will now use the description of (pseudo-) Riemannian
almost-product structures and  almost-complex structures
with a Norden metric in terms of a pair of metrics (twin or dual metrics).
Let us consider a triple $(M,h,K)$ where $M$ is a differentiable manifold,
$h$ is a metric on $M$ and $K$ is a (1,1) tensorfield on $M$ such that
the following holds
$$K^{2}=\epsilon I,  \ \
K^t hK=\ep h $$
As we said in the Introduction such a triple defines
a  (pseudo-)Riemannian
almost-product structure if $\ep =1$, while it defines an almost-complex
structure with a Norden metric $h$ if $\ep =-1$. The triple $(M,h,K)$ admits
an equivalent description as another triple $(M,h,g)$, where
$g$ is a metric on $M$ satisfying the relation
$(h^{-1}g)^2 =\ep I $ or the equivalent relation
$(g^{-1}h)^2 =\ep I $, because of the following elementary
proposition.\smallskip\\

{\bf Proposition 2.2.}\ {\em   Let $K$ and $h$
be real  $n\times n\ $ matrices and $\ep=+1$ or $-1$. Then the matrices $K$
and $h$ satisfy the  relations}
$$h^t=h,\ \ det\  h\neq 0,\ \ K^{2}=\epsilon I,\ \  K^t hK=\ep h \eqno	(2.20)
$$
{\em if and only if there exists a real matrix $g$ such that $g$ and $h$
satisfy  the relations}
$$
h^t=h,\ \  det\  h\neq 0, \ \ g^t=g, \ \  (h^{-1}g)^2 =\ep I \eqno  (2.21)
$$
{\em Moreover one has
$$
g=hK \eqno	    (2.22)
$$
and}
$$
( K^{-1})^{2}=\epsilon I,\ \  K^{-1t} gK^{-1}=\ep g  \eqno	    (2.23)
 $$

{\it Proof}. If one has (2.20) then define $g$ by (2.22) and check (2.21)
and (2.23). Conversely, if one has (2.21) then define $K=h^{-1}g$ and
check (2.20). The claim is proved. \hfill (Q.E.D.)\smallskip\\

We can therefore state the following theorem.\smallskip\\

{\bf Theorem 2.3.}\ {\em  Let $M$ be a $n$-dimensional manifold, $n>2$,
with a metric $h$ and a symmetric connection $\Ga$ and let us consider
the Euler-Lagrange equations (2.3)-(2.4) for the action (2.1). Let us assume
that the analytic function $f(S)$ is such that eq. (2.5)
has an isolated root $S=c$ with $f^{\prime}(c) \neq 0$;
setting then $g_{\mu\nu}=R_{(\mu\nu)}(\Ga)$, the Euler-Lagrange equations
imply the  relations
$(h^{-1}g)^2 =\ep I $,\ $Ric(g)=g $. Therefore:
\begin{itemize}
\item [(i)] if $\ep>0$ after rescaling and denoting
$P=g^{-1}h $ one gets an almost-product Einstein manifold $(M,g,P),$ i.e.
 $$
 Ric(g)=\ga g
$$
$$
P^2=I,\ \ g(PX,PY)=g(X,Y),\ \ X,Y\in \chi (M)\ \ ;$$
\item [(ii)] if  $\ep<0$ after rescaling and denoting $J=g^{-1}h $
 one gets instead an  anti-Hermitian  Einstein manifold $(M,g,J)$, i.e.
 $$
 Ric(g)=\ga g
$$
$$
J^2=-I,\ \ g(JX,JY)=-g(X,Y),\ \ X,Y\in \chi (M)
$$\end{itemize}
Here $\chi (M)$ is the Lie algebra of vector fields on $M$.}\smallskip\\

Notice	that the signatures of the metrics $h$ and $g$ in the almost-product
case can be in principle
arbitrary, while in the almost-complex case the signature is $(m,m)$.
In any case they will be lower--semicontinuous functions, so that  without
any restriction we can
assume that they remain constant in the neighborhood of a
generic point; in particular they will not change in connected components
of the manifold. In the next section we will therefore
study the equation $(h^{-1}g)^2 =\ep I$ for a generic point on the manifold,
i.e. study it algebraically as a matrix equation.

{\it Remark.} Notice that $\ep=0$ corresponds to the case  $S=0$, which holds
iff $f(0)=0$. Then we
have  to distinguish  two subcases $f^\prime (0)=0$ and $f^\prime (0)\neq 0$.
In the first subcase
both equations (2.3) and (2.4) are automatically satisfied and no condition
for $g$ and $\Ga$ arise. Therefore any pair $(g, \Ga)$ solves this subcase.
When  $f^\prime (0)\neq 0$ then (2.3) is instead equivalent to the algebraic
equation $[h^{-1} S(\Ga)]^2=0$, which leads to an {\it almost-tangent}
structure \cite{Eli}, while (2.4) remains unchanged. We shall not discuss
this case in the present paper.

\section{Solutions of the Matrix Equations}
\ \ \
\indent Let us then consider the matrix equation
$$
(\zz)^2 = \ep \uno_n \eqno (3.1)
$$
where $h$ and $g$ are symmetric non-degenerate
real $n\times n\ $ matrices, $\uno_n$ is the
identity matrix in $n$ dimensions and $\ep$ is a
non-vanishing real number.

In order to solve equation (3.1) we first notice that
it is manifestly  ${\rm Gl}(n,\R)$--invariant under the canonical right-action
$
(h,g) \mapsto (A^th A, A^tgA)
$
where $A^t$ denotes the matrix transpose to $A$.
More exactly, transforming $h$ and $g$ as  metrics
one can observe that $P=\zz$ transforms by a similarity
transformation $P \rightarrow A^{-1} P A$
(i.e. as a (1,1) tensor).
Equation (3.1) is also invariant under the transformation
$
(h, g, \ep)\mapsto (g, h, \ep^{-1})
$.
Moreover, we have
$
(\det(g)/\det(h))^2 = \ep^n
$, so that when
$n$ is even there are no restrictions on the sign
of $\ep$, while $\ep$ has to be positive when $n$ is odd.

It is always possible to rescale $g$ by $\sqrt{\mid\ep\mid}$ and reduce
(3.1) to the canonical form $(\zz)^2 = \pm \uno_n$.

When $\ep$ is positive equation (3.1) admits always the trivial solution
$g=\sqrt\ep\, h$; however, this does not exhaust all the possible solutions.
Let us first observe in fact, by standard
minimal-polynomial arguments (see e.g. \cite{HJ}), that the  matrix equation
$$P^2=\uno_n \eqno (3.2)$$
admits only  solutions of the form $P=M^{-1} D_k M$ for some non-singular
matrix $M$, where the matrices $D_k$ (Jordan forms) are diagonal
$$
D_k=
\left( \ba{cc}
-\uno_k & 0 \\ 0 & \uno_{n-k} \ea \right) \eqno (3.3)
$$
and $k=0,\ldots, n$. This result can be restated as follows:
{\it if an automorphism $P$ of the
vector space $\R ^n$ satisfies (3.2) then there exists a basis
in which $P$ is represented by one	of the matrices $D_k$}. The
non-negative integer $k$ is an
invariant of this authomorphism. In fact such an automorphism $P$ ($P^2=id$)
represents an
almost--product structure on $\R ^n$. The set of all solutions of the equation
$(\zz)^2 = \uno_n$ is then described by the following
theorem (compare \cite{HJ} Th. 4.5.15 case II b):\smallskip\\

{\bf Theorem 3.1.}\ {\em
 Let $g=g^t$ and $h=h^t$ be two real (symmetric) non-degenerate
matrices (metrics). Then the
following are equivalent:
\begin{description}
\item[a)] $(\zz)^2 = \uno_n$
\item[b)] the two metrics $h$ and $g$ are simultaneously diagonalizable
with $\pm 1$ on the diagonal, i.e.
there exists a real non-degenerate matrix $R$ such that
$$
h=R^t D_h R,\ \ \ \ \ g=R^t D_g R
$$
and $D_h$ and $D_g$ are diagonal matrices with $+1$ or $-1$ on the diagonal.
\end{description}}
The proof of this theorem is given in the Appendix.\smallskip\\

Let us proceed to discuss the case $(\zz)^2=-\uno_n$. It is known that if
$J$ is any $n\times n$ real matrix satisfying the relation
$$
J^2=-\uno_n \eqno (3.4)
$$
then $n$  must be even, $n=2m$, $J$ can be represented	as
$$J=M\J M^{-1}$$
where $\J$ is the canonical form
$$
\J=
\left( \ba{cc}
0 & \uno_m \\ -\uno_m & 0 \ea \right) \eqno (3.5)
$$
and $M$ is a nondegenerate real matrix. In fact one deals with a complex
structure and the matrix $\J$ gives the canonical complex structure on
$\R ^{2m}$. The following holds true: \smallskip\\

{\bf Theorem 3.2.}\ {\em
Let $h=h^t$ and $g=g^t$ be two $2m\times 2m$ real (symmetric) non-degenerate
matrices (metrics). Then the  following are equivalent:
\begin{description}
\item[a)] $(\zz)^2 = -\uno_{2m}$
\item[b)] there exists a real non-degenerate matrix $R$ such that
$$
h=R^t \left( \ba{cc}
\uno_m & 0 \\ 0 & -\uno_m \ea \right)  R,\ \ \ \ \ g=R^t \left( \ba{cc}
0 & \uno_m \\ \uno_m & 0 \ea \right)  R
$$
i.e. in the appropriate coordinate system the two metrics $g$ and $h$
take the following canonical forms
$$
K_h= \left( \ba{cc}
\uno_m & 0 \\ 0 & -\uno_m \ea \right)  ,\ \ \ \ \ \ K_g= \left( \ba{cc}
0 & \uno_m \\ \uno_m & 0 \ea \right) .
$$
\end{description} }
Also this proof is given in the Appendix.\smallskip\\

From theorems 3.1 and 3.2 it follows that locally for any given
metric $g$ one can construct a twin metric $h$. If $g$ satisfies
$Ric(g)=g$ this means that locally one produces an almost-product or an
almost-complex Einsteinian structure.

\section{$K$-Structures and K\"ahler-like Manifolds}
\ \ \
\indent We first present here a formalism which at once describes properties
of various structures important in differential geometry, such as
 almost-complex  and almost-product structures, Hermitian and
 anti-Hermitian metrics, K\"ahler manifolds and locally decomposable
 manifolds, {\it etc...}, and then, in the next Sections, consider in more
detail the pseudo-K\"ahlerian  and anti-K\"aherian metrics.

 Let $M$ be a smooth manifold, $TM$ its tangent bundle
and $\chi (M)$ the algebra of vectorfields on $M$.
A {\it $(K,\ep)$-structure}  ({\it $K$-structure} in short) on $M$ is a field
of endomorphisms $K$ on
$TM$ such that $K^2=\ep I$, where $\ep =\pm 1$. Thus $\ep=1$ corresponds to
an almost-product
structure, while $\ep=-1$ provides an almost-complex structure.\\

Let $\nabla : \chi (M)\times \chi (M)\rightarrow \chi (M)$
be  a connection, denoted by $(X,Y)\rightarrow \nabla_XY$.
A $K$-structure is {\it integrable} iff there exists
an a linear torsionless connection on $M$ such that $\nabla K=0$;
 or equivalently  the Nijenhuis tensor $N$
$$
N(X,Y)=[KX,KY]-K[KX,Y]-K[X,KY]+\ep [X,Y]
$$
vanishes. If $K$ is integrable then there exists an atlas of adapted
coordinate charts
on $M$ in which $K$ takes a canonical form (see {\it e.g.} \cite{Kob}).\\

{\bf Definition 1.}\ A 5-tuple $(M,K,g,\ep,\si)$ is called a
{\it $(K,g)$-manifold} if $g$ is a metric on $M$ and $K$ is a $K$-structure,
$K^2=\ep I$, such that
$$
g(KX,KY)=\si g(X,Y)\eqno (4.1)
$$
for all vectorfields $X$ and $Y$ on $M$. Here $\si =\pm 1$. In this case we
shall say that the metric $g$ is {\it K-compatible } (or a {\it K-metric} in
short).\smallskip\\
The definition above unifies the following four cases:
the case $\ep=1,\si=1$ corresponds to a (pseudo-) Riemannian almost-product
structure;
the case $\ep=1,\si=-1$ provides an almost para-Hermitian structure; the case
$\ep=-1,\si=1$ is  known as an almost-(pseudo)-Hermitian structure; and
finally the case $\ep=-1,\si=-1$  corresponds to an almost-complex
structure with a Norden metric.\\

Introduce a (0,2) tensorfield $h$, the {\it twin} of $g$, by
$$
h(X,Y)=g(KX,Y) \eqno (4.2)$$
Then
$$
h(X,Y)=\ep \si h(Y,X),\ \ \ \ \ h(KX,KY)=\si h(X,Y) \eqno (4.3)$$
Notice that for $\ep \si=1$ the twin tensor is a metric (and this is, in fact,
the case we have obtained from our variational principle), while for
$\ep \si=-1$ the twin tensor is a two-form (and one deals with an
almost-Hermitian or almost para-Hermitian structure).

Let $\psi$ be a $(0,3)$ tensorfield defined by the formulae
$$
\psi (X,Y,Z) = g((\na_X K)Y,Z)\equiv (\na_X h)(Y,Z)  \eqno (4.4)$$
In a coordinate language $\psi_{\al\mu\nu}$ is nothing but
$\na_\al h_{\mu\nu}$. It possesses the following properties
$$
\psi (X,Y,Z) = -\si \psi (X,KY,KZ) = \si \ep \psi (X,Z,Y) \eqno (4.5)$$
Notice that the classification of almost-Hermitian structures \cite{GH},
Riemannian almost-product structures \cite{Nav} as well as almost-complex
structures with a Norden metric \cite{GB} is based on algebraic properties
of $\psi$: namely, one decomposes $\psi$ into irreducible components under
the action of the appropriate group.
The most restrictive class is K\"ahler-like, when simply $\psi=0$.
If the tensorfield $\psi$ vanishes then automatically  $\na K=0$ for a
torsionless (Levi-Civita)
connection and the Nijenhuis tensor $N$ is forced to vanish, too
(\cite{Yan,YK} c.f. also formulae
 (6.7)).  Therefore the corresponding $K$-structure is integrable. It leads
to the following definition:\smallskip\\

{\bf Definition 2.}\  A metric $g$ on a $(K,g)$-manifold is called a
{\it K\"ahler-like metric} if $\nabla K=0$, i.e.
$$
\nabla_{X}(KY)=K\nabla_XY,\ \  X,Y\in\chi (M)\eqno(4.6)
$$
where $\nabla$ is the Levi-Civita connection of $g$ itself.\smallskip\\

If $\ep =-1, \si =1$ then a $K$-metric is called a K\"ahler metric.
If $\ep=1, \si=1$ then a $K$-metric shall be called
a {\it pseudo-K\"ahlerian metric} (it is also called  a (pseudo-) Riemannian
locally decomposable metric). The case $\ep=1,\si=-1$ shall be considered in a
forthcoming publication \cite{BFV} (see also \cite{RMG,BBR}). Our results on
anti-K\"ahlerian manifolds ($\ep=-1, \si=-1$) will be presented in sections 6
 and 7. The following Proposition extends the Proposition 3.6 in
\cite{KN} to an arbitrary $(K,g,\ep,\si)$-structure\smallskip\\

{\bf Proposition 4.1.}\ {\em The Riemann curvature $R(X,Y)Z$ and the Ricci
tensor $S(X,Y)$ of a K\"ahler-like manifold (M,K,g) satisfy the
following properties:
$$
R(X,Y)\circ K=K\circ R(X,Y) ,\ \ \ \ R(KX,KY)=\si R(X,Y) \eqno (4.7)$$
$$
S(KX,KY)=\si S(X,Y),\ \ \ (\si-\ep)S(X,Y) =\hbox{tr}[V\mapsto K(R(X,KY)V)]
\eqno (4.8) $$
}
{\it Proof:}\ The proof is a simple repetition of the proof of
Proposition 3.6 in \cite{KN}, provided one suitably takes into account
formulae (4.3). \hfill (Q.E.D.)\smallskip\\

Consider now the twin $F$ of the Ricci tensor $S$
$$
F(X,Y)=S(KX,Y) \eqno (4.9)$$
Then
$$
F(X,Y)=\ep \si F(Y,X),\ \ \ F(KX,KY)=\si F(X,Y) \eqno (4.10)$$
Notice that the symmetry property of $F$ is exactly the same as for $h$.
Therefore, we conclude the following: \smallskip\\

{\bf Lemma 4.2.}\ {\em A K\"ahler-like manifold is Einstein iff the dual of
$S$ is proportional to the dual of $g$, {\it i.e.} $F(X,Y)\sim h(X,Y)$.\\
}
\smallskip\\
Notice that for $\ep\si=-1$ both twins $F$ and $h$ are two-forms.
This Lemma in the K\"ahlerian case leads to a necessary condition on the
first Chern class for a manifold to have an Einstein-K\"ahler metric
\cite{Cal,Yau}. Recall the Goldberg conjecture \cite{Gol}
(see also \cite{Ols,SV}) saying that almost-K\"ahler Einsteinian manifold is
a complex one. It means that an Einstein almost-Hermitian manifold with a
closed K\"ahler form is automatically Hermitian, {\it i.e.} its
almost-complex structure is integrable. An extension of the Goldberg
conjecture to the other K\"ahler-like manifolds will be disscuss
elsewhere \cite{BFV}.

\section{Almost Product Einstein Manifolds}
\ \ \
\indent In this section we consider the problems of the global
existence and classification of  almost--product  Einstein  manifolds.
At the	begining we shall recall some basic facts about almost-product and
(pseudo-) Riemannian almost-product structures.

 The simplest examples of almost-product structures
are product manifolds, {\it i.e.} manifolds which are the Cartesian
product of two manifolds
$$
M = M_1 \times M_2 \eqno (5.1)$$
In this case the tangent bundle splits as $TM = TM_1\oplus TM_2$ and
$P=P_2 - P_1$,
where $P_i$ are the corresponding projections on $TM_i$, $i=1,2$.
More generally, giving
an almost-product structure is equivalent to splitting	the tangent bundle
into two complementary subbundles ({\it distributions} or
{\it almost-foliations}):
$TM = V\oplus H$; in this case $P= P_V - P_H$. $P$ is integrable iff
both the distributions are integrable. ({\it i.e.} they are {\it foliations}).
Integrable almost-product structures are also called  {\it locally
product manifolds} \cite{Yan} since locally they have the form (5.1). It
means that locally (around any point), there exists an adapted coordinate
system $(x^a, y^\al)$, $a=1,...,k$ and $\al=1,...,n-k$, such that the
tensorfield $P$ takes the canonical form (3.3), {\it i.e.}
$P\pt_a=-\pt_a$ and $P\pt_\al=\pt_\al$.

Similarly, one can consider other structures: for example, a (pseudo-)
Riemannian product
manifold as a product of two (pseudo-) Riemannian manifolds
$$
(M, g) = (M_1, g_1)\times (M_2, g_2) \eqno (5.2)$$
where $g=g_1\oplus g_2$.
More generally, an almost-product (pseudo-) Riemannian structure
$(M, P, g)$ is integrable iff $\na^g P = 0$ for the Levi--Civita
connection $\na^g$ of $g$.
In this case we speak of a {\it locally decomposable} (pseudo-)Riemannian
manifold; in the present note we shall however propose to call it a
{\it pseudo-K\"ahler manifold}. For locally decomposable (pseudo-)
Riemannian structures both  foliations
are {\it totally geodesic} \cite{Yan,YK}. In an adapted coordinate system
the metric $g$ "separates the variables" (see \cite{Yan,YK})
$$
ds^2=g_{ab}(x)dx^a dx^b + g_{\al\bt}(y) dy^\al dy^\bt \eqno (5.3)$$
The twin metric has the form $h=h_1\ominus h_2$, {\it i.e.}
$$
ds_h^2=g_{ab}(x)dx^a dx^b - g_{\al\bt}(y) dy^\al dy^\bt $$
Since $\na^g_\mu h_{\al\bt}\equiv \psi_{\mu\al\bt}=0$ then we have in this
case $\Ga(g)=\Ga(h)$. In fact, this property gives an equivalent definition
of pseudo-K\"ahler manifolds, provided $g$ and $h$ are twin metrics.

Recall that an almost-product Einstein manifold is a triple $(M,g,P)$ where
$g$ is a metric and $P$ is a $(1,1)$ tensorfield which satisfy
$$
 Ric(g)=\ga g	  \eqno (5.4)
$$
$$
P^2=I,\ \ g(PX,PY)=g(X,Y),\ \ X,Y\in \chi (M) \eqno (5.5)
$$
Notice that $P=I$ and $P=-I$ give trivial examples of an almost-product
structure. Non-trivial examples are given by the following:\smallskip\\

 {\bf Proposition 5.1.}\
{\em Let ($M^n,g)$ be an Einstein manifold satisfying (5.4) with an
indefinite metric $g$ of signature $q$, $1\leq q <n$ . Then  there exists
on $M^n$ a non-trivial almost-product structure $P$ satisfying (5.5) and
therefore one gets an almost-product Einstein manifold $(M,g,P)$.}

{\it Proof:}\  If $g$ is a pseudo-Riemannian metric on $M$
then it was proved in \cite{CM}
that there exist a (strictly) Riemannian metric $h$ and an almost product
structure $P$ on $M$ such that $g(X,Y)=h(PX,Y)$ and $h(PX,PY)=h(X,Y)$
for all vectorfields $X$ and $Y$ on $M$. Therefore one gets the relations
(5.5). The almost-product structure $P$ is non-trivial since if $P=\pm I$
then $g=\pm h$; but in fact the metric $g$ is pseudo-Riemannian while $h$
is strictly Riemannian.\hfill (Q.E.D.)\smallskip\\

It follows from the proposition above that any	manifold with
a strictly pseudo-Riemannian
Einstein metric serves as an example of a pseudo-Rieman\-nian
almost-product Einsteinian manifold.

It should also be noted that construction of $P$ for a given
pseudo-Riemannian metric $g$ is not a canonical one (and, in fact,
depends on a choice of some "background" Riemannian
metric on $M$). Therefore, a single (possibly Einstein) pseudo-Riemannian
metric leads, in principle,
to several almost-product (Einsteinian) manifolds. It makes a striking
difference between the solutions of our variational problem corresponding
to the positive roots of the fundamental
equation (2.5) and those corresponding to the negative roots. In the second
case, as we shall see in the next Section, there are further topological
obstructions for the existence of an almost-complex structure.

Let M be a pseudo-K\"ahler manifold, {\it i.e.} (locally) in adapted
coordinate systems $(x^a, y^\al)$ the metric $g$ splits as
$g=g_1\oplus g_2$, where $g_{ab}=g_1(x)$ and
$g_{\al\bt}=g_2(y)$. If both metrics in (5.3) are Einstein
$$R_{ab}(g_1)=\ga g_{ab}, \ \ \ R_{\al\bt}(g_2)=\ga g_{\al\bt}\eqno (5.6)$$
for the same constant $\ga$, then it follows:
$$R_{AB}(g) = \ga g_{AB}  \eqno (5.7)$$
Therefore, one has\smallskip\\

{\bf Proposition 5.2.}\  {\em A pseudo-K\"ahlerian manifold  is Einstein
iff in any adapted coordinates $(x^a,y^\al)$ both metrics are Einstein
for the same constant $\ga$.
}

{\it Proof:} See \cite{Yan,YK}.
\hfill (Q.E.D.)\smallskip\\

Interesting examples of locally product (pseudo-) Riemannian manifolds
which are not locally decomposable are given by {\it warped product}
spacetimes \cite{One,DVV,CC}. Given two (pseudo-) Riemannian manifolds
$(M_i, g_i)$, $i=1,2$, and a smooth function
$\theta :M_1\longrightarrow \R$, on the product manifold $M=M_1\times M_2$
put the metric $g=g_1\oplus e^{2\theta} g_2$. The resulting (pseudo-)
Riemannian manifold $M=M_1\times_\theta M_2$ is called a
{\it warped product} manifold. It is of course an almost-product
(pseudo-) Riemannian manifold and it is
conformal to a locally decomposable one. It is an interesting and
intriguing fact that many exact solutions of Einstein equations
(including {\it e.g.}  Schwarzschild,
Robertson-Walker, Reissner-Nordstr\"om, de Sitter {\it etc...}) and
also p-brane solutions (\cite{AVV}) are, in fact, warped product
spacetimes \cite{CC}. Therefore, these exact solutions provide beautiful
examples of almost-product Einsteinian manifolds.

There are topological restrictions on a (paracompact) manifold $M$ for
the existence of an almost-product structure  of rank $k$, which
are the same as for the existence of
a (strictly) pseudo-Riemannian metric of signature $(k, n-k)$, which are
again the same as for the existence of a $k$-dimensional distribution.
For example, for the existence of a metric with Lorentz signature on a
compact manifold $M$ ({\it i.e.} for the existence of a nowhere vanishing
vectorfield)  the necessary and sufficient
condition is that the Euler characteristic number vanishes.

\section{Anti-K\"ahlerian Manifolds}
\ \  \
\indent Now we consider in some detail the case
of a $K$-metric with
 $\ep=-1, \si=-1$, which we call an anti-K\"ahlerian metric.\smallskip\\

 {\bf Definition 3.}\ A triple $(M,g,J)$, where $J$ is an almost-complex
structure and the metric $g$ is anti-Hermitian:
$g(JX,JY)=-g(X,Y),X,Y\in \chi(M)$
 is called an {\it anti-K\"ahlerian manifold} if $\nabla J=0$, where
 $\nabla$ is the Levi-Civita connection.\smallskip\\

 We will prove that a real anti-Hermitian metric on a complex manifold
satisfies the K\"ahler condition $\nabla J=0$
on the same manifold treated as a real manifold if and only if the metric
is the real part of a holomorphic metric on this manifold.

 Let $(M,J)$ be a  $2m$-dimensional almost-complex  real manifold
and let $g$ be an anti-Hermitian  metric on $M$. We extend $J$,
$g$  and the Levi-Civita connection $\nabla$ by $\C$-linearity
to the complexification of the tangent bundle $T_{C} M=T_M\otimes \C$.
We use the same notation for the complex extended $g, J$ and $\na$.
Then the Levi-Civita connection is the mapping $(X,Y)\rightarrow\nabla_XY$
where $X$ and $Y$ are now complex vectorfields	({\it i.e.}, sections of
$T_{C} M$). Then the (complex extended) torsion tensor $T$ vanishes
$$
T(X,Y)=\nabla_XY-\nabla_YX-[X,Y]=0
$$
and the ordinary formulae are valid for the connection:
$$
\nabla_Xg(Y,Z)=Xg(Y,Z)-g(\nabla_XY,Z)-g(Y,\nabla_XZ)=0\eqno (6.1)
$$
and for the Riemann tensor
$$
R(X,Y)Z=[\nabla_X,\nabla_Y]Z-\nabla_{[X,Y]}Z\eqno (6.2)
$$
where $X,Y,Z$ are  complex vectorfields. For the sake of clarity, we
stress that for the moment we are just complexifying the tangent bundle
but  we do not assume the almost-complex structure $J$
is  integrable. Let us now fix a (real) basis
$\{X_1,...,X_m,JX_1,...,JX_m\} $ in each tangent space
$T_xM$; then the set $\{Z_a, Z_{\bar a}\}$, where
$Z_a=X_a-iJX_a,\  Z_{\bar a}=X_a+iJX_a$,
forms a basis for each complexified tangent space $T_x M\otimes\C$.
Unless otherwise stated, little Latin indices $a,b,c,...$ run from
$1$ to $m$, while Latin capitals $A,B,C,...$
run through
$1,...,m,{\bar 1},...,{\bar m}$; for notational convenience we shall also
bar capital indices and we
shall assume ${\bar{\bar A}}=A$. One has $JZ_a=iZ_a$
and $J Z_{\bar a}=-i Z_{\bar a}.$
We set
$
g_{AB}=g( Z_A,Z_B)=g_{BA}$. Then the following holds:\smallskip\\

{\bf Proposition 6.1.}\ {\em Let $(M, J)$ be an almost-complex manifold
and $g$ be an anti-Hermitian
metric on it. Then the complex extended  metric $g$ (in the complex basis
constructed above)
satisfies the following conditions
$$
g_{a\bar b}=g_{\bar ba}=0\eqno (6.3)
$$
$$
g_{\bar A\bar B}={\bar g_{AB}}\eqno (6.4)
$$}

{\it Proof.}\
Since the metric $g$ is anti-Hermitian, we have
$$
g(Z_a, Z_{\bar b})=-g(JZ_a,J Z_{\bar b})=
-g(iZ_a,- i Z_{\bar b})=-g( Z_a,  Z_{\bar b})
$$
Therefore  $g_{a\bar b}=0$, which proves (6.3). The proof of (6.4)
is well known. In fact we have
$$
g_{\bar a\bar b}=g(Z_{\bar a},Z_{\bar b})=g(X_a+iJX_a,X_b+iJX_b)
$$
$$
=g(X_a,X_b)-g(JX_a,JX_b)+ig(JX_a,X_b)+ig(X_a,JX_b)
$$
and
$$
g_{ab}=g(X_a-iJX_a,X_b-iJX_b)
$$
$$=
g(X_a,X_b)-g(JX_a,JX_b)-ig(JX_a,X_b)-ig(X_a,JX_b)
$$
Therefore we get $g_{\bar a\bar b}={\bar g}_{ab}$.
Similarly we consider the other components of $g_{ AB}$ and hence we prove
(6.4).\hfill (Q.E.D.)\smallskip\\

 It is customary to write a metric satisfying (6.3) and (6.4) as
$$
ds^2=g_{ab}dz^adz^b+g_{\bar a\bar b}dz^{\bar a}dz^{\bar b}  \eqno (6.5)
$$

We define now the complex Christoffel symbols $\Ga^C_{AB}$ as
$$
\nabla_{Z_A}Z_B=\Ga^C_{AB}Z_C\eqno (6.6)
$$
It is known \cite{Yan,KN} that if $\nabla J=0$ then the torsion $T$ and
the Nijenhuis tensor $N$ satisfy the identity
$$
T(JX,JY)=\frac{1}{2}N(X,Y) \eqno (6.7)
$$for any vectorfields $X$ and $Y$. Since  the complex extended
Levi-Civita connection $\nabla$ has no torsion, the complex Christoffel
symbols are symmetric.
  In this case the complex structure $J$ is integrable so that the real
manifold $M$ inherits the structure of a complex manifold. Let us now
recall (see {\it e.g.} \cite{KN}) that there is a
one-to-one correspondence between complex manifolds and real manifolds
with an integrable complex structure.
 This means that there exist real, adapted (local) coordinates
$(x^1,...,x^m,y^1,...,y^m)$ such that
 $$
 J(\frac{\partial}{\partial x^a})=\frac{\partial}{\partial y^a},\ \
   J(\frac{\partial}{\partial y^a})=-\frac{\partial}{\partial x^a}
$$
Setting $z^a=x^a+iy^a$ and taking $X_a=\partial/\partial x^a$ one gets
$$
Z_a=X_a-iJX_a=2\partial/\partial z^a=2\partial_a,\ \  Z_{\bar a}=X_a+iJX_a
=2\partial/\partial {\bar z}^a=2\partial_{\bar a}
$$
where $\partial_A=\partial/\partial z^A$ and $z^{\bar a}={\bar z}^a$. It
appears that
$z^a$'s form a complex (analytic) coordinate chart on $M$.
Now from (6.1) one gets
$$
\Ga^C_{AB}=\frac{1}{2}g^{CD}(Z_Ag_{DB}+Z_Bg_{DA}-Z_Dg_{AB})
=g^{CD}(\partial_Ag_{DB}+\partial_Bg_{DA}-\partial_Dg_{AB})\eqno (6.8)
$$

 Notice that the  relation (6.1) is valid for the complex
 extended metric $g$ and complex vector fields $X,Y,Z$ if and only if
 it is valid for real vectorfields.\smallskip\\

{\bf Theorem 6.2.}\ {\em Let $M$ be a $m$-dimensional complex manifold,
thought as a real $2m$-dimensional manifold with a complex structure $J$.
Let us further assume that $M$ is provided with an anti-Hermitian
metric	$g$. We extend	 $J$, $g$  and the Levi-Civita connection
$\nabla$ by $\C$-linearity
to the complexified tangent bundle $T_{C}M$. Then the following conditions
are equivalent:

$(i)$
$$
\nabla_{X}(JY)=J\nabla_XY     \eqno(6.9)
$$
where $X$ and $Y$ are arbitrary real vectorfields;

$(ii)$
the (complex) Christoffel symbols satisfy
$$
\Ga ^C_{AB}=0\ \ \ \hbox{except} \ \hbox{for}\ \ \Ga_{ab}^c\
\hbox{ and}\ \ \Ga_{\bar
a\bar b}^{\bar c }={\bar \Ga_{ab}^c}	\eqno (6.10)
$$

$(iii)$
there exists  a local complex coordinate system
$(z^1,...,z^m)$  on $M$ such that the components of the complex
extended metric  $g_{ab}$ in the canonical  form (6.5) are holomorphic
functions
$$
\partial_{{\bar c}}\, g_{ab}=0	 \eqno (6.11)
$$}

{\it Proof.}\  From (6.6)  we have
$${\bar \Ga}^C_{AB}=
\Ga^{\bar C}_{{\bar A}{\bar B}}
$$
The connection satisfies  the conditions
$$\nabla_{Z_B}(JZ_c)=J\nabla_{Z_B}Z_c=i\nabla_{Z_B}(Z_c)$$
$$\nabla_{Z_B}(JZ_{\bar c})=J\nabla_{Z_B}Z_{\bar c}=
-i\nabla_{Z_B}(Z_{\bar c})$$
if and only if
$$
\Ga^{ a}_{B{\bar c}}=\Ga^{\bar a}_{Bc}=0 \eqno (6.12)
$$
This proves the equivalence between $(i)$ and $(ii)$.
Then for the Christoffel symbols (6.8) by taking  (6.3) into account
one gets
$$
\Ga^{ a}_{b{\bar c}}=g^{aD}(\partial_bg_{D\bar c}+\partial_{\bar c}g_{Db}
-\partial_Dg_{D\bar c})=
g^{ad}\partial_{\bar c}g_{bd} \eqno (6.13)
$$
and from (6.12) it follows that
$$
\partial_{\bar c}g_{bd}=0 \eqno (6.14)
$$

The other relations (6.12)  also are reduced to (6.14) or its complex
conjugated. Therefore the relation (6.14) is equivalent to (6.11).
This proves the equivalence  between $(i)$
 and  $(iii)$. Our claim is thence proved. \hfill (Q.E.D.)\smallskip\\

We have proved that a real anti-Hermitian metric on a complex manifold
satisfies the K\"ahler condition $\nabla J=0$
on the same manifold treated as a real manifold if and only if the metric
is the real part of a holomorphic metric on this manifold. Therefore,
there exists a one-to-one correspondence between anti-K\"ahler manifolds
and complex Riemannian manifolds with a holomorphic metric as they were
defined in \cite{Bru} (see also \cite{GI}).

From an algebraic viewpoint, let us mention that we have been
dealing  with the following construction. Let $V$ be a real vector
space with a complex structure $J$ and let $G$ be a complex-valued
bilinear form on $V$. Let us set
$$
F(X,Y)=G(X,Y)-G(JX,JY)-iG(X,JY)-iG(JX,Y)
$$
Then we have
$$
F(JX,JY)=-F(X,Y)
$$
Now one takes the real	(or imaginary) part of $F$ to get
a real anti-Hermitian bilinear form on $V$.


\section{ Anti K\"ahlerian Einstein manifolds}
\ \ \
\indent In this section we consider the problems of the global
existence and classification of  anti-Hermitian  Einstein  manifolds.
Recall that an anti-Hermitian Einstein manifold is a triple $(M,g,J)$
where $g$ is a metric and $J$ is a $(1,1)$ tensorfield which satisfy
 $$
 Ric(g)=\ga g	  \eqno (7.1)
$$
$$
J^2=-I,\ \ g(JX,JY)=-g(X,Y),\ \ X,Y\in \chi (M) \eqno (7.2)
$$
Then the metric $g$ has necessarily  the signature $(m,m)$ (see section 3),
being $2m=\hbox{dim}M$. Let us show that by taking the real part of a
holomorphic Einstein metric  on a complex manifold of  complex dimension
$m$ one can get a real Einstein  manifold of real dimension $2m$.

From (4.8) we have for the Ricci tensor
$$
Ric(g)(JX,JY)=-Ric(g)(X,Y) \ \ X,Y\in\chi(M)
$$
Therefore, analogously to (6.3), we have
$$
R_{a\bar b}=0 \eqno (7.3)
$$

We shall not attempt here to consider solutions of Einstein equations for
a generic metric of the form (6.5) but consider only the case when
$g_{ab}$ is a holomorphic function
$$
\partial_{\bar c}\,g_{ab}=0  \eqno (7.4)
$$
From (7.4) and (4.7) we get for the Riemann tensor:
$$
R^D_{ABC}=0\ \ \ \hbox{except} \ \hbox{for}\ \ R^d_{abc}\ \hbox{ and}\ \
R^{\bar d}_{\bar a\bar b\bar c}={\bar R^d_{abc}} \eqno (7.5)
$$
The (complex) Einstein equations
$$
R_{AB}(g)=\ga g_{AB} \eqno (7.6)
$$
are thence equivalent to a pair of equations
$$
R_{ab}(g_{cd})=\ga g_{ab}  \eqno (7.7a)
$$
$$
R_{\bar a\bar b}(g_{\bar c\bar d})=\ga g_{\bar a\bar b}  \eqno (7.7b)
$$
To get a real solution of Einstein equations (7.1) from (7.7)
one uses real coordinates $(x^{\mu}), \mu =1,...,2m$ on $M,$ i.e.
 $z^a=x^a+ix^{m+a}, a=1,...,m$ and writes the metric (6.5) as
$$
ds^2=g_{ab}dz^adz^b+g_{\bar a\bar b}dz^{\bar a}dz^{\bar b}
=g_{\mu\nu}dx^{\mu}dx^{\nu}		   \eqno (7.8)
$$
where $g_{\mu\nu}$  is a real metric. We have thence proved the following
theorem:\smallskip\\

{\bf Theorem 7.1.}\ {\em If $(M,g,J)$ is an anti-K\"ahlerian manifold,
i.e. a complex manifold of complex
dimension $m$ with a holomorphic metric $g_{ab}(z), a,b=1,...,m$ and a
real metric $g_{\mu\nu}(x),
\mu,\nu=1,...,2m$ defined by (7.8), then the holomorphic metric  $g_{ab}(z)$
satisfies (7.7a) if and only if the real metric $g_{\mu\nu}(x)$ is a
solution of the Einstein equations (7.1)
$$
R_{\mu\nu}(g)=\ga g_{\mu\nu}	\eqno (7.9)
$$
}

As an example one can take a complex analytic continuation of any
real analytic solution of Einstein equations. A simple example is
$$
ds^2=dz^adz^a+\frac{(z^adz^a)^2}{1-z^az^a}+\ \hbox{
complex conj.}=g_{\mu\nu}dx^{\mu}dx^{\nu} \eqno (7.10)
$$
This metric $ g_{\mu\nu}$ on "the complex sphere": $w_1^2+...+w_{m+1}^2=1$
(which can be interpreted as a quadric
$\zeta_1^2+...+\zeta^2_{m+1}-\zeta_{m+2}^2=0$
in $\C P^{m+1}$ if one takes $w_i=\zeta_i/\zeta_{m+2,}$)
gives a solution of the Einstein equations (7.8) and provides an example
of an anti-Hermitian Einstein manifold $(M,g,J).$

In particular for $m=2$ we get a real solution of Einstein equations
on  the 4-dimensional real manifold ($w_1^2+w_2^2+w_3^2=1, w_i\in \C$)
with a metric of signature $(++--)$.

Notice also that any Einstein metric
on a compact Riemannian manifold $M^n$ leads to an anti-K\"ahlerian
Einstein metric on another real manifold ${\cal M}^{2n}$. It follows
from the known fact \cite{Bou} that any Einstein metric is analytic
in a certain atlas on $M^n$. Therefore there exists a complex analytic
continuation of the metric to a complex manifold of complex dimension
$n$ which is a real anti-K\"ahlerian manifold  ${\cal M}^{2n}$.

\section*{Acknowledgments}
\ \ \
\indent Two of us (A.B. and I.V.) gratefully acknowledge
the hospitality of the Institute of Mathematical Physics ``J.--L. Lagrange''
of the University of Torino and
the support of G.N.F.M. of
Italian C.N.R. This work is sponsored by G.N.F.M., M.U.R.S.T. (40\% \ Proj.
``Metodi Geometrici e Probabilistici in Fisica Matematica''); one of us
(A.B.) acknowledges also the support from KBN 2 P302 023 07.

\section*{Appendix}
\ \  \
\indent

{\it Proof of Theorem 3.1}:\  The proof $b)\Rightarrow a)$ follows obviously
from the fact that each two diagonal matrices with $\pm 1$ on their
diagonals do satisfy our equation, which is
invariant under the appropriate transformation.\\
The converse $a)\Rightarrow b)$ is less obvious. As we already know
there exists a real non-degenerate matrix $M$ such
$$\zz=M D_k M^{-1} \eqno (A.1) $$
From this one gets
$$
\g=\h D_k
$$
where $\g=M^t gM$ and $\h=M^t hM$. Since $\g$, $\h$ and $D_k$ are symmetric
matrices  one has thence
$$
\h D_k=D_k\h  \eqno (A.2)
$$
Let us now represent $\h$ in block--form:
$$
\left( \ba{cc}
\h_{11} & \h_{12}\\ \h^t_{12} & \h_{22} \ea \right)
$$
where $\h^t_{11}=\h_{11}$, \ \,$\h^t_{22}=\h_{22}$ and $\h_{11}$ is a
$k\times k$ matrix. Then from (A.2) we obtain\\
$$
\left( \ba{cc}
\h_{11} & -\h_{12}\\ \h^t_{12} & -\h_{22} \ea \right)=
\left( \ba{cc}
\h_{11} & \h_{12}\\ -\h^t_{12} & -\h_{22} \ea \right)
$$
Therefore $\h_{12}=0$ and one gets
$$
\h=\left( \ba{cc}
\h_{11} & 0\\ 0 & \h_{22} \ea \right),\ \ \ \
\g=D_k\h=\left( \ba{cc}
\h_{11} & 0\\ 0 & -\h_{22} \ea \right)	\eqno (A.3)
$$
Now we make use of the fact that any real nondegenerate
symmetric matrix is ${}^t$-congruent
to a diagonal matrix whose diagonal elements are equal to $+1$ or $-1$, i.e.:
$h_{11}=S^t_1 D_{h_{11}}S_1$ \
(and analogously $h_{22}=S^t_2 D_{h_{22}}S_2$), \,where
$D_{h_{11}}$ (resp. $D_{h_{22}}$)  is a diagonal matrix
with $\pm 1$ along the diagonal.
Therefore one has
$$
\h=S^t\left( \ba{cc}
D_{h_{11}} & 0\\ 0 & D_{h_{22}} \ea \right) S,\ \ \ \
\g=S^t\left( \ba{cc}
D_{h_{11}} & 0\\ 0 & -D_{h_{22}} \ea \right) S
$$
where $S=\left( \ba{cc}
S_1 & 0 \\ 0 & S_2 \ea \right) $. Notice that $S$
commutes with $D_k$, \,i.e. $SD_k=D_k S$.\\
Taking then $R=SM$ the theorem is proved.
\medskip\hfill (Q.E.D.)\smallskip\\

{\it Proof of Theorem 3.2}:\ To prove $b)\Rightarrow a)$ check that
$K^{-1}_h K_g =\J$ and
then make use of the appropriate transformation properties.

In order to prove the converse $a)\Rightarrow b)$, in full analogy
with the previous case one makes use of the fact that $\zz=M\J M^{-1}$.
Now this leads to the condition (with the same
notation) :
$$
\h\J=-\J\h \eqno (A.4)
$$
Writing then $\h$ in block--form:
$$
\h=\left( \ba{cc}
a & b \\ b^t & d \ea \right)
$$
one has from (A.4):
$$
\left( \ba{cc}
-b & a \\ -d & b^t \ea \right)=\left( \ba{cc}
b^t & d \\ -a & -b \ea \right)
$$
Therefore $a=-d$, $b=b^t$ and we have
$$
\h=\left( \ba{cc}
a & b \\ b & -a \ea \right) \eqno (A.5a)
$$
$$
\g=\h\J=\h=\left( \ba{cc}
-b & a \\ a & b \ea \right) \eqno (A.5b)
$$
Let us further show that there exists a real matrix $S$ such that
$$
S\J=\J S \eqno (A.6)
$$
$$
\h=S^t K_h S \eqno (A.7a)$$
$$\g=S^t K_g S \eqno (A.7b)
$$
To be sure that (A.6) holds take the matrix $S$ in the form
$$
S=\left( \ba{cc}
s & u \\ -u & s \ea \right)
$$
where $s$ and $u$ are real $m\times m$ matrices to be determined from the
conditions (A.7a) and (A.7b). Equation (A.7a) reads as
$$ \left( \ba{cc}
s^ts-u^tu & s^tu+u^ts \\ u^ts+s^tu & u^tu-s^ts \ea \right)=\left( \ba{cc}
a & b \\ b & -a \ea \right)    \eqno (A.8a)
$$
while the equation (A.7b) gives
$$
\left( \ba{cc}
-u^ts-s^tu & s^ts-u^tu \\ s^ts-u^tu & s^tu+u^ts \ea \right)=\left( \ba{cc}
-b & a \\ a & b \ea \right)   \eqno (A.8b)
$$
Notice that eq. (A.8a) is equivalent to (A.8b) so that we are left with the
following conditions
$$
s^ts-u^tu=a,\ \ \ \ s^tu+u^ts=b
$$
This last equation can be rewritten in the complex form
$$
(s+iu)^t (s+iu) = a+ib
$$
Now it is known (see for example \cite {HJ}) that any nondegenerate
symmetric complex matrix $a+ib$ can be represented in the form
$$
a+ib=N^tN
$$
where $N$ is a (nondegenerate) complex	matrix. Taking $s+iu=N$ and
$R=SM^{-1}$ the theorem is proved.\hfill (Q.E.D.)\smallskip\\

\end{document}